\documentclass[11pt,a4paper,reqno]{article}
\usepackage{graphicx} 
\usepackage{authblk}
\usepackage{fullpage}
\usepackage[margin=1in]{geometry}
\usepackage{parskip}
\usepackage{hyperref}
\usepackage{eqnarray,amsthm, amssymb, amsmath,amsfonts,mathrsfs,verbatim,epsfig}
\usepackage{enumerate}
\usepackage{tikz} 
\usepackage{pgfplots}
\usetikzlibrary{calc,matrix,trees,arrows.meta}

\newtheorem{theorem}{Theorem}[section]
\newtheorem{lemma}[theorem]{Lemma}
\newtheorem{corollary}[theorem]{Corollary}

\newtheorem{remark}{Remark}[section]
\newtheorem{conjecture}{Conjecture}[section]

\newtheorem{proposition}[theorem]{Proposition}
\newtheorem{problem}[theorem]{Problem}

\newcommand{\fref}[1]{Figure \ref{fig:#1}}

\newcommand{\cref}[1]{Conjecture \ref{conjecture:#1}}

\newcommand{\N}{\mathbb{N}}

\newcommand{\al}{\alpha}
\newcommand{\la}{\lambda}

\newcommand{\sg}{\sigma}
\newcommand{\ga}{\gamma}

\newcommand{\SCP}{\mathrm{SCP}}
\newcommand{\sh}{\mathrm{sh}}
\newcommand{\cont}{\mathrm{cont}}
\newcommand{\wt}{\mathrm{ht}}

\newcommand{\Sn}[1]{\mathcal{S}_{#1}}

\newcommand{\T}{\mathcal{T}}

\newcommand{\inc}{{\mathrm{inc}}}

\newcommand{\latticemn}[3]{
	\foreach \x in {1,...,#1} {
		\draw[line width=1pt] ($(#3)+(-\x+1,-\x+1)$) -- +(#2-1,-#2+1);
		\foreach \y in {1,...,#2} {
			\draw[fill] ($(#3)+(\y-\x,-\x-\y+2)$) circle (3pt);
		}
	}
	\foreach \y in {1,...,#2} {
		\draw[line width=1pt] ($(#3)+(\y-1,-\y+1)$) -- +(-#1+1,-#1+1);
	}
}

\newcommand{\latticemnud}[3]{
	\foreach \x in {1,...,#1} {
		\draw[line width=1pt] ($(#3)+(\x-1,\x-1)$) -- +(-#2+1,#2-1);
		\foreach \y in {1,...,#2} {
			\draw[fill] ($(#3)+(-\y+\x,\x+\y-2)$) circle (3pt);
		}
	}
	\foreach \y in {1,...,#2} {
		\draw[line width=1pt] ($(#3)+(-\y+1,\y-1)$) -- +(#1-1,#1-1);
	}
}

\newcommand{\latticemnlabel}[3]{
	\foreach \x in {1,...,#1} {
		\foreach \y in {1,...,#2} {
			\node at ($(#3)+(-\y+\x-#1+#2,\x+\y-#1-#2-2+1.25)$) {\scriptsize $(\x,\y)$};
		}
	}
}
\newcommand{\latticemnlabelrank}[4]{
	\foreach \x in {1,...,#1} {
		\foreach \y in {1,...,#2} {
			\ifnumequal{#4}{\x+\y-2}{
				\node at ($(#3)+(-\y+\x-#1+#2,\x+\y-#1-#2-2+1.1)$) {\scriptsize $(\x,\hspace{-0.5mm}\y)$};
				\filldraw ($(#3)+(-\y+\x-#1+#2,\x+\y-#1-#2)$) circle (6pt);
				
			}{}
		}
	}
}

\newcommand{\latticemnp}[4]{
	\foreach \x in {1,...,#1} {
		\foreach \y in {1,...,#2} {
			\draw[line width=1pt] ($(#4)+(\y-\x,-\x-\y+1)$) -- +(0,1-#3);
		}
	}
	\foreach \z in {1,...,#3} {
		\latticemn{#1}{#2}{$(#4)+(0,-\z)$}
	}
}

\newcommand{\ptnF}[2]{
	\foreach \x [count=\i] in {#1} {
		\draw[line width=1pt] (#2)+(0,\i-1) grid +(\x,\i);
	}
}

\newcommand{\tabbound}[2]{
	\coordinate (start) at (#2);
	\foreach \d [count=\i] in {#1} {
		\ifnum\i=1
		\draw[line width=1.5pt, line cap=rect] (#2) -- ++(\d,0) -- ++(0,1) coordinate (current);
		\xdef\te{\d}
		\else
		\draw[line width=1.5pt, line cap=rect] (current) -- ++(\d-\te,0) -- ++(0,1) coordinate (current);
		\xdef\te{\d}
		\fi
	}
	\draw[line width=1.5pt, line cap=rect] (current) -- ++(-\te,0) -- (#2);
}
\newcommand{\specialribbonh}[2]{
	\draw[line width=1pt, rounded corners] ($(#1)+(-0.5,-0.5)$) circle (2pt) \foreach \d [count=\i] in {#2} {\ifodd\i {-- ++($(\d,0)$)} \else {-- ++(0,-\d)} \fi} circle (2pt);
}
\newcommand{\specialribbonv}[2]{
	\draw[line width=1pt, rounded corners] ($(#1)+(-0.5,-0.5)$) circle (2pt) \foreach \d [count=\i] in {#2} {\ifodd\i {-- ++($(0,-\d)$)} \else {-- ++(\d,0)} \fi} circle (2pt);
}

\newcommand{\colorchain}[3]{
	\draw[rounded corners,line cap=round,#3,line width=10pt,opacity=0.5] (#1) \foreach \d [count=\i] in {#2} {\ifodd\i {-- ++($(\d,-\d)$)} \else {-- ++(-\d,-\d)} \fi};
}
\newcommand{\colorchainl}[3]{
	\draw[rounded corners,line cap=round,#3,line width=10pt,opacity=0.5] (#1) \foreach \d [count=\i] in {#2} {\ifodd\i {-- ++($(-\d,-\d)$)} \else {-- ++(\d,-\d)} \fi};
}
\newcommand{\colorchainlnop}[3]{
	\draw[rounded corners,line cap=round,#3,line width=10pt,nearly opaque] (#1) \foreach \d [count=\i] in {#2} {\ifodd\i {-- ++($(-\d,-\d)$)} \else {-- ++(\d,-\d)} \fi};
}

\newcommand{\colorchainld}[3]{
	\draw[rounded corners,line cap=round,#3,line width=10pt,opacity=0.5] (#1) \foreach \d [count=\i] in {#2} {\ifodd\i {-- ++($(-\d,-\d)$)} \else {-- ++(0,-\d)} \fi};
}

\newcommand{\PP}[1]{
	\foreach \x [count=\i] in {#1}{
        \ifnum\i>1
        \times
        \fi
        \mathbf{\x}
	}
}

\title{Stanley's conjecture on the Schur positivity of distributive lattices}

\author[1]{Grace M.X. Li
	\thanks{\texttt{grace\_li@sust.edu.cn}}
}
\author[2]{Dun Qiu
	\thanks{\texttt{qiudun@nankai.edu.cn}}
}
\author[3]{Arthur L.B. Yang
	\thanks{\texttt{yang@nankai.edu.cn}}
}
\author[4]{Zhong-Xue Zhang
	\thanks{\texttt{zhzhx@mail.nankai.edu.cn}}
}
\affil[1]{School of Mathematics and Data Science, 
	Shaanxi University of Science and Technology, Xi’an, Shaanxi 710021, P. R. China}
\affil[2,3,4]{Center for Combinatorics, LPMC, Nankai University, Tianjin 300071, P. R. China}

\date{}

\begin{document}
\maketitle

\noindent\textbf{Abstract.} 
In this paper we solve an open problem on distributive lattices, which was proposed by Stanley in 1998. This problem was motivated by a conjecture due to Griggs, which equivalently states that the incomparability graph of the boolean algebra $B_n$ is nice.   
Stanley introduced the idea of studying the nice property of a graph by investigating the Schur positivity of its corresponding chromatic symmetric functions.  Since the boolean algebras form a special class of distributive lattices, Stanley raised the question of  whether the incomparability graph of any distributive lattice is Schur positive.  Stanley further noted that this seems quite unlikely.  
In this paper, we construct a family of distributive lattices which are not nice and hence not Schur positive. We also provide a family of distributive lattices which are nice but not Schur positive.

\noindent\textbf{AMS Classification 2020: } 05E05, 06A07 


\section{Introduction}
The main objective of this paper is to study the Schur positivity of the chromatic symmetric functions of the incomparability graphs of distributive lattices. Stanley stated in his 1998 paper \cite{Sta98} that perhaps  every distributive lattice is Schur positive, which implies that every distributive lattice is nice.  In this paper we provide several families of distributive lattices whose chromatic symmetric functions are not Schur positive.  Some of our examples showed that distributive lattices need not be nice.  This solves Stanley's open problem in the negative.   

Let us first give an overview of related backgrounds.  Let $B_n$ denote the boolean algebra of order $n$, which is the poset of all subsets of $[n]=\{1,2,\ldots,n\}$, ordered by inclusion.  The celebrated Sperner's theorem states that $B_n$ has width $\binom{n}{\lfloor \frac{n}{2}\rfloor}$, namely, the maximal size of antichains is equal to the size of its middle rank.   One can associate a partition $\lambda$ of $2^n$ to each chain decomposition of $B_n$ by rearranging the sizes of chains in weakly decreasing order,  in which case we also say that this chain decomposition is of type $\lambda$.  Given two partitions $\mu$ and $\nu$ of the same size, we say that $\mu$ is less than or equal to $\nu$  in dominance order, denoted by $\mu\unlhd  \nu$, if $\sum_{i=1}^k \mu_i\leq \sum_{i=1}^k \nu_i$  for any $k\geq 1$.  Griggs \cite{Griggs} proposed the following conjecture. 
\begin{conjecture}[{\cite[Problem 3]{Griggs}}]\label{conj-griggs}
    If  a partition $\lambda$ of $2^n$ is less than or equal to the partition corresponding to the symmetric chain decomposition of $B_n$  in dominance order, then there is a chain decomposition of type $\lambda$.  
\end{conjecture}
Stanley \cite{Sta98} proposed an approach to the above conjecture by studying the Schur positivity of the chromatic symmetric function of the incomparability graph of the boolean algebra $B_n$.  
For a given poset $P$, its incomparability graph $\inc(P)$ has the vertex set $P$ and the edge set  $\{(i,j):i,j\textnormal{ are incomparable in }P\}$.  In his seminal paper \cite{Sta95}, Stanley introduced the concept of the chromatic symmetric function of a graph. For a given graph $G$ with vertex set $V$ and edge set $E$,  a proper coloring $\kappa$ of $G$ is a map from $V$ to $\{1,2,\ldots\}$ such that $\kappa(u)\neq \kappa(v)$ if $(u,v)\in E$.  The chromatic symmetric function of $G$ is defined by
$$
X_G=X_G({\bf x})=\sum_{\kappa} \prod_{v\in V} x_{\kappa(v)},
$$
where ${\bf x} = \{x_1,x_2,\ldots\}$ and $\kappa$ ranges over all proper colorings of $G$.   
When $x_1=\cdots=x_n=1$ and $x_i=0$ for $i>n$, the symmetric function $X_G$ gives the chromatic polynomial $\chi_G(n)$ of Birkhoff \cite{Birkhoff12}.  

Stanley \cite{Sta95} studied the expansions of $X_G$ in terms of various bases of symmetric functions. 
In this paper we are interested in the expansion of  $X_G$ in terms of Schur function basis  $\{s_\la\}$. Note that the Schur functions are of special importance in the theory of symmetric functions, due to their strong connection with irreducible representations of  symmetric groups. See Fulton and Harris \cite{Fulton}, Stanley \cite{EC2} or Sagan \cite{Sag} for further details. Recall that if a symmetric function $f$ can be written as a nonnegative linear combination of certain basis $\{b_\la\}$ then it is said to be $b$-positive.  The elementary symmetric function basis $\{e_\la\}$ and the complete homogeneous symmetric function basis $\{h_\la\}$ are both $s$-positive bases related to some classical representations.  By convention, we also say that a graph $G$ (or a poset $P$) is Schur positive or $s$-positive if the chromatic symmetric function $X_G$ (or $X_{\inc(P)}$) is $s$-positive.

There has been a substantial amount of research about the $s$-positivity and further, $e$-positivity of chromatic symmetric functions. One of the most famous problems in combinatorics is the $(3+1)$-free conjecture of Stanley and Stembridge \cite{SS93}, which states that the chromatic symmetric function of the incomparability graph of any $(3+1)$-free poset is $e$-positive, and it has been open for more than thirty years. Haiman \cite{Hai} and Gasharov \cite{Ga96} proved that these symmetric functions are $s$-positive. Gasharov first raised the Schur positivity conjecture of claw-free graphs in an unpublished work, see Stanley \cite{Sta98}, which also remains open. 

Stanley \cite{Sta98} proposed a promising conjecture on the Schur positivity of the boolean algebra $B_n$, which  is stronger than Conjecture \ref{conj-griggs} in the sense that we shall describe below. For a graph $G$, a \textit{stable partition} of $G$ is a partition of its vertex set so that any two vertices in the same part are disconnected.  The type of a stable partition is defined to be the partition $\la$ obtained by rearranging the sizes of parts in weakly decreasing order.  
Following Stanley \cite{Sta98}, we say that a graph  $G$ is \textit{nice} if $G$ has a stable partition of type $\la$, then it also has a stable partition of type $\mu$ for each $\mu\unlhd\la$.  We also say that a poset $P$  is nice if so is $\inc(P)$.  Thus Conjecture \ref{conj-griggs} equivalently states that the boolean algebra $B_n$ is nice for any $n$. Stanley established the following connection between the Schur positivity and the nice property for graphs. 
\begin{proposition}[{\cite[Proposition 1.5]{Sta98}}] \label{prop-spos-nice}
    If $G$ is Schur positive, then it is nice. 
\end{proposition}
Motivated by Conjecture \ref{conj-griggs} and the above proposition,  
Stanley \cite{Sta98} proposed the following conjecture.
\begin{conjecture}[{\cite[p.\ 270]{Sta98}}]\label{conjecture:StanleyBn} For any $n\in\N$, the incomparability graph
    $\inc(B_n)$ is Schur positive.
\end{conjecture}
Since the boolean algebras $B_n$ are special distributive lattices, Stanley \cite{Sta98} further raised the following problem. 
\begin{problem}[Stanley, 1998]\label{conjecture:Stanley}
     Is $X_{\inc(L)}$ Schur positive for any distributive lattice  $L$?
\end{problem}
Let $\hat{L}$  be the poset obtained from a distributive lattice $L$ by removing the maximal element $\hat{1}$ and the minimal element $\hat{0}$. Stanley gave a distributive lattice $L$ which is Schur positive  but $\hat{L}$ is not Schur positive.   Based on his example, Stanley  remarked that the Schur positivity of $X_{\inc(L)}$ might not be true for all distributive lattices  $L$.   

In this paper we give a negative answer to Problem \ref{conjecture:Stanley}.  The remainder of this paper is organized as follows. In Section \ref{section-2}, we include some definitions and results on distributive lattices and chromatic symmetric functions which will be used in subsequent sections.  In Section \ref{section-3} we construct a family of distributive lattices which are not nice and hence not Schur positive by Proposition \ref{prop-spos-nice}. Lonc and Elzobi \cite{LE99} proved that the product of two chains is always a nice distributive lattice.  In Section \ref{section-4} we prove that the product of two chains are not Schur positive in general. 
In Section \ref{section-6}, we summarize a few open problems for further research.





\section{Preliminaries} \label{section-2}
The aim of this section is two-folds. Firstly, we present some basic definitions and results on distributive lattices, where we adopt the terminologies and notations of Davey and Priestley \cite{DB}.
Secondly, we follow Wang and Wang \cite{WW20} to review the Schur function expansion of the chromatic symmetric functions in terms of special rim hook tabloids, which will be frequently used later.  

\subsection{Distributive lattices} 
A \textit{lattice} is a poset $L$ in which for any pair of elements $a,b\in L$, the join $a\vee b$ and the meet $a\wedge b$ exist.  A subposet  $M\subseteq L$ is a \textit{sublattice} of $L$ if  for any $a,b \in M$, we have $a\vee b \in M$ and $a\wedge b\in M$.  A lattice $L$ is said to be a \textit{distributive lattice} if the operations $\vee$ and $\wedge$ satisfy that for any $a,b,c \in L$,
\begin{eqnarray*}
    a\wedge(b\vee c) &=& (a\wedge b)\vee(a\wedge c),\\
    a\vee(b\wedge c) &=& (a\vee b)\wedge(a\vee c).
\end{eqnarray*}
For any positive integers $n,n_1,n_2,\ldots,n_r$, let $\mathbf{n}$ be the $n$-element chain $1\leq 2\leq\cdots\leq n$, 
and let $\mathbf{n_1}\times \mathbf{n_2}\times \cdots\times \mathbf{n_r}$ be the product of $r$ chains.  
It is routine to check that $\mathbf{n_1}\times \mathbf{n_2}\times \cdots\times \mathbf{n_r}$ is a graded distributive lattice of rank $n_1+n_2+\cdots+n_r -r$.  In particular, when $n_i=2$ for all $i=1,\ldots,r$, $\mathbf{n_1}\times \mathbf{n_2}\times \cdots\times \mathbf{n_r}$ gives the boolean algebra $B_r$. 
\fref{posets} gives examples of posets $\mathbf{6}$, $\mathbf{6}\times \mathbf{3}$ and $\mathbf{6}\times \mathbf{2}\times \mathbf{2}$. 
\begin{figure}[ht!]
    \centering
    \begin{tikzpicture}[scale=0.5]
        \latticemn{6}{1}{-7,-1};
        \latticemn{6}{3}{0,0};
        \latticemnp{6}{2}{2}{9,1}
    \end{tikzpicture}
    \caption{Posets $\mathbf{6}$, $\mathbf{6}\times \mathbf{3}$ and $\mathbf{6}\times \mathbf{2}\times \mathbf{2}$.}\label{fig:posets}
\end{figure}
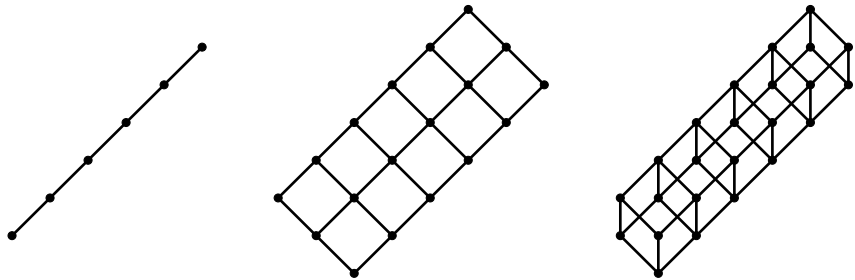

 
    

As mentioned before, the nice property of a distributive lattice is mainly concerned with its chain decompositions, also known as its chain partitions.  Recall that a  \textit{chain partition} of a poset $P$ is a partition $\{C_1, \ldots,C_k\}$ of the set $P$, so that the induced subposet $P[C_i]$ is a chain for each $i=1,\ldots,k$.  In this paper we also use the notion of {semi-ordered chain partitions}.  By a \textit{semi-ordered chain partition} we mean a chain partition such that all parts are ordered according to their sizes in decreasing order and the parts of the same size are also ordered. Let $\SCP_{P,\la}$ be the set of semi-ordered chain partitions of $P$ of type $\la$. For example, $\SCP_{\PP{4},(2,1,1)}$ is
\begin{align*}
    \{&(\{1,2\},\{3\},\{4\}),\quad (\{1,2\},\{4\},\{3\}),\quad (\{1,3\},\{2\},\{4\}),\quad
    (\{1,3\},\{4\},\{2\}),\\
    &(\{1,4\},\{2\},\{3\}),\quad (\{1,4\},\{3\},\{2\}),\quad
    (\{2,3\},\{1\},\{4\}),\quad (\{2,3\},\{4\},\{1\}),\\
    & (\{2,4\},\{1\},\{3\}),\quad (\{2,4\},\{3\},\{1\}),\quad 
    (\{3,4\},\{1\},\{2\}),\quad (\{3,4\},\{2\},\{1\})\}.
\end{align*}

\subsection{Schur expansion of chromatic symmetric functions}

We assume that the readers are familiar with the theory of symmetric functions. Here we only present some necessary definitions and results which will be used in the Schur expansion of chromatic symmetric functions. We refer the readers to  Macdonald \cite{Macd-book}, Mendes  and Remmel \cite{MR15}, and Stanley \cite{EC2} for more information on symmetric functions. 

Let $\mathbb{N}$ be the set of natural numbers. By a \textit{weak composition} of $n\in \mathbb{N}$ we mean a sequence $\beta=(\beta_1,\,\ldots,\,\beta_\ell)\in \mathbb{N}^{\ell}$ such that  $\beta_1+\cdots+\beta_\ell=n$, denoted by $\beta\vDash n$.
A \textit{partition} $\la$ of $n$ is a weakly decreasing sequence $\la=(\la_1,\,\ldots,\,\la_\ell)$ of positive integers such that $\la_1+\cdots +\la_\ell=n$. Let $\la^{\langle k \rangle}$ denote the partition obtained from $\la$ by removing its first $k-1$ parts, namely, 
$\la^{\langle k \rangle}=(\la_{k},\,\ldots,\,\la_\ell)$.
For any partition $\la$ of $n$ (denoted by $\la \vdash n$), we also write $\lambda=\langle 1^{\alpha_1},\, 2^{\alpha_2},\, \ldots,\, n^{\alpha_n}\rangle$, where $\alpha_k$ is the number of parts of size $k$. 
Given any partition $\la=(\la_1,\ldots,\la_\ell)\vdash n$, the \textit{Ferrers diagram} (or Young diagram) 
of $\la$ is the graph of $n$ boxes arranged in $\ell$ left-justified rows, such that there are $\la_i$ boxes in row $i$ (indexed from bottom to top). A \textit{rim hook} (or ribbon) of a partition $\la$ is a set of connected boxes of its Ferrers diagram containing no $2\times2$ shape, such that its removal also gives a Ferrers diagram of a partition. If one continues removing rim hooks from the remaining shape until no boxes are left, then these removed rim hooks form a partition of the original Ferrers diagram, called a \textit{rim hook tabloid} of shape $\la$. A \textit{special rim hook tabloid} is a rim hook tabloid such that each rim hook has a cell in the first column. See \fref{rimhook} for examples of shape $(5,3,2,1)$.

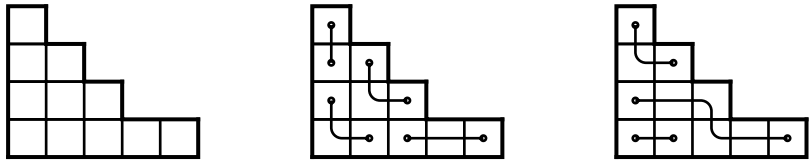
\begin{figure}[ht!]
    \centering
    \begin{tikzpicture}[scale=0.5]
        \ptnF{5,3,2,1}{0,0}
        \tabbound{5,3,2,1}{0,0}
        
        \ptnF{5,3,2,1}{8,0}
        \tabbound{5,3,2,1}{8,0}
        \specialribbonv{9,4}{1}
        \specialribbonv{10,3}{1,1}
        \specialribbonv{9,2}{1,1}
        \specialribbonh{11,1}{2}
        
        \ptnF{5,3,2,1}{16,0}
        \tabbound{5,3,2,1}{16,0}
        \specialribbonv{17,4}{1,1}
        \specialribbonh{17,2}{2,1,2}
        \specialribbonh{17,1}{1}
    \end{tikzpicture}
    \caption{A Ferrers diagram, a rim hook tabloid and a special rim hook tabloid.}\label{fig:rimhook}
\end{figure}

For a (special) rim hook tabloid  $T$,  its shape is denoted by $\sh(T)$. The \textit{content} of $T$ is the partition of $n$ formed by sizes of rim hooks in $T$, denoted by $\cont(T)$. Notice that if $T$ is a special rim hook tabloid then $\sh(T)\unlhd\cont(T)$ in dominance order \cite{ER}.
Given partitions $\la,\mu\vdash n$, let $\T_\la$ be the set of special rim hook tabloids of shape $\la$, and let $\T_{\la,\mu}$ be the subset of $\T_\la$ consisting of those with content $\mu$. It is clear that $\T_{\la,\mu}\neq\emptyset$ only if $\la\unlhd\mu$. 
For each rim hook $B$ of some $T\in\T_\la$, its height $\mathrm{ht}(B)$ is defined to be the number of rows it spans minus one, and the height $\mathrm{ht}(T)$ is defined to be the sum $\sum_{B\in T} \mathrm{ht}(B)$. 
For example, the special rim hook tabloid $T$ in \fref{rimhook} has $\wt(T)=2$ and $\cont(T)=(6,3,2)$. 

Let $[f]_g$ be the coefficient of $f$ in $g$, and let $\{m_\la\}$ be the monomial symmetric function basis. As a consequence of the Jacobi-Trudi formula, E\v{g}ecio\v{g}lu and Remmel \cite{ER} proved the following result. 
\begin{theorem}[{\cite[Theorem 1]{ER}}]\label{theorem:ER} 
    For any $\la,\mu\vdash n$, we have
    $$[s_\la]_{m_\mu} = \sum_{T\in\T_{\la,\mu}} (-1)^{\wt(T)}.$$
\end{theorem}
According to the definition of chromatic symmetric functions, Wang and Wang \cite{WW20} generalized the result of E\v{g}ecio\v{g}lu and Remmel and gave a Schur expansion formula for chromatic symmetric functions. To state their result, we also need the notion of semi-ordered stable partitions of graphs. Similar to the definition of semi-ordered chain partitions of posets, by a \textit{semi-ordered stable partition} of a graph $G$ we mean a {stable partition} such that all parts are arranged according to their sizes in weakly decreasing order and the parts of the same size are also ordered.
The following result gives a combinatorial interpretation of the coefficients in the Schur expansion of chromatic symmetric functions. 

\begin{theorem}[{\cite[Theorem 3.1]{WW20}}]\label{theorem:WW}
    For any graph $G=(V,E)$ and any partition $\la\vdash |V|$, we have
    $$
    [s_\la]_{X_G} = \sum_{T\in\T_\la} (-1)^{\wt(T)} N_T,
    $$
    where $N_T$ is the number of semi-ordered stable partitions of G of type $cont(T)$.
\end{theorem}

For example, taking graph $G=(V,E)$ with $V=\{1,2,3\}$ and $E=\emptyset$, we have 
\begin{align}\label{eq-s-expansion}
X_G = e_1^3 = s_{3}+2s_{2 1}+s_{1 1 1}.
\end{align}
Note that there are two special rim hook tabloids $T_1,T_2\in\T_{(2,1)}$ as shown in \fref{sr21}.
One can check that $\wt(T_1)=0$,\, $N_{T_1} = |\{(\{1,2\},\{3\}),(\{1,3\},\{2\}),(\{2,3\},\{1\})\}|=3$, $\wt(T_2)=1$ and $N_{T_2} = |\{(\{1,2,3\})\}|=1$. 
Thus, by Theorem \ref{theorem:WW}, we have $[s_{2 1}]_{X_G}=3-1=2$, which coincides with the above expansion \eqref{eq-s-expansion}.

\begin{figure}[ht!]
    \centering
    \begin{tikzpicture}[scale=0.5]
        \ptnF{2,1}{0,0}
        \tabbound{2,1}{0,0}
        \specialribbonv{1,2}{}
        \specialribbonh{1,1}{1}
        \node at (-1,1) {$T_1$};
        
        \ptnF{2,1}{6,0}
        \tabbound{2,1}{6,0}
        \specialribbonv{7,2}{1,1}
        \node at (5,1) {$T_2$};
    \end{tikzpicture}
    \caption{Two special rim hook tabloids in $\T_{(2,1)}$.}\label{fig:sr21}
\end{figure}
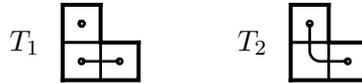

Note that the semi-ordered chain partitions of a poset $P$ are in bijection with the semi-ordered stable partitions of its incomparability graph $\inc(P)$. 
Thus, for the incomparability graph of a poset $P$, Theorem \ref{theorem:WW} can be rephrased as follows.  

\begin{corollary}\label{coro:WW}
For any poset $P$ with $n$ elements and any partition $\la\vdash n$,
\begin{equation}\label{formula}
[s_\la]_{X_{\inc(P)}} = \sum_{T\in\T_\la} (-1)^{\wt(T)} |\SCP_{P,\cont(T)}|.
\end{equation}
\end{corollary}

\section{Not all distributive lattices are nice}  \label{section-3}

In this section, we construct a family of distributive lattices which are not nice. 
By Proposition \ref{prop-spos-nice}, these lattices are not either Schur positive. 

For any $n\geq 1$, let $B_{3,n}$ be the poset on the set $\{a,\ b,\ c,\ d,\ e,\ f,\ 1',\ 1,\ 2',\ 2,\ \ldots,\ n',\ n\}$ with the Hasse diagram as shown in \fref{B_{3,n}}. 
Note that $B_{3,1}$ is just the boolean algebra of order $3$. 

\begin{figure}[ht]
    \begin{center}
    \begin{tikzpicture}[scale=0.5]
        \draw[line width=1pt] (0,0) --(1,1)--(2,0)--(1,-1)-- cycle;
        \draw[line width=1pt] (0,1) --(1,2)--(2,1)--(1,0)-- cycle;
        \draw[line width=1pt] (0,0)--(0,1);
        \draw[line width=1pt] (1,1)--(1,2);
        \draw[line width=1pt] (2,0)--(2,1);
        \draw[line width=1pt] (1,-1)--(1,0);
      \draw[fill](1,1) circle(2pt);
      \draw[fill](2,0) circle(2pt);
      \draw[fill](1,-1) circle(2pt);
      \draw[fill](0,1) circle(2pt);
      \draw[fill](1,2) circle(2pt);
      \draw[fill](2,1) circle(2pt);
      \draw[fill](1,0) circle(2pt);
      \node at (1,2.3) [font=\small]{$a$};
      \node at (2.3,1) [font=\small]{$d$};
      \node at (0.7,1) [font=\small]{$c$};
      \node at (-0.5,1) [font=\small]{$b$};
      \node at (2.3,0) [font=\small]{$f$};
      \node at (0.7,-0.1) [font=\small]{$e$};
       \node at (-0.5,0) [font=\small]{$1$};
       \node at (-1.5,-0.75) [font=\small]{$2$};
       \node[rotate=135] at (-2.4,-1.7) [font=\small]{$\vdots$};
       \node at (-3.6,-2.6) [font=\small]{$n-1$};
       \node at (-4.4,-3.7) [font=\small]{$n$};
       \node at (1.5,-1.25) [font=\small]{$1'$};
       \node at (0.5,-2.25) [font=\small]{$2'$};
       \node[rotate=135] at (-0.7,-3.25) [font=\small]{$\vdots$};
       \node at (-1.0,-4.25) [font=\small]{$n-1'$};
       \node at (-2.5,-5.25) [font=\small]{$n'$};
        \latticemn{5}{2}{0,0};
    \end{tikzpicture}
    \end{center}
    \caption{$B_{3,n}$.}\label{fig:B_{3,n}}
    \end{figure}
    
First we have the following result.
\begin{lemma}
For any positive integer $n$, the poset $B_{3,n}$ is a distributive lattice.
\end{lemma}

\begin{proof}
It is known that any sublattice of a distributive lattice is still distributive \cite[page 293]{EC1}.
One can check that $B_{3,n}$ is a sublattice of the distributive lattice  $(\mathbf{n+1})\times \mathbf{2}\times\mathbf{2}$, hence it is distributive.  
\end{proof}

The main result of this section is the following.

\begin{theorem}\label{main-thm-not-nice}
 For each $n\geq 6$ the distributive lattice $B_{3,n}$ is not nice.
\end{theorem}
\begin{proof}
Let $\la=(n+3,\, n+1,\, 2)$ and let $\mu \vdash 2n+6$ be a partition of length $3$ such that $\mu_1=n$. 
Note that if $n\geq 6$ then such a partition $\mu$ always exists, say $\mu=(n,n,6)$. It is also clear that  $\mu\unlhd  \la$ in dominance order.  Observe that  $\{\{a,\ d,\ f,\ 1',\ \cdots,\ n'\},\ \{c,\ 1,\ \cdots,\ n\},\ \{b,\ e\}\}$ is a chain partition of $B_{3,n}$ of type $\la$.   We proceed to show that there exists no chain partition of $B_{3,n}$ of type $\mu$. 
Assume to the contrary there exists a chain partition $\{C_1,\ C_2,\ C_3\}$ of type $\mu$.
Since $\{f,\ e,\ 1\}$ is a $3$-element antichain, $|\{f,\ e,\ 1\}\cap C_i|=1$ for $i=1,2,3$. Without loss of generality, assume that $1\in C_1$ and hence $i\in C_1$ for any $2\le i\le n$ since each such $i$ is incomparable with $f$ and $e$. 
Since $\{d,\ b,\ c\}$ is a $3$-element antichain, either $b$ or $c$ must belong to $C_1$ since $d$ and $1$ are incomparable. Hence $|C_1|\geq n+1$, contradicting to the assumption that $\mu_1\le n$. 
\end{proof}

\begin{remark}
When $n=1$, $B_{3,1}=B_3$ is the boolean algebra verified to be Schur positive. We have also checked using Sagemath that $B_{3,n}$ is nice whenever $n\leq 5$.
\end{remark}
Now we see that distributive lattices do not necessarily have nice property. However, due to the work of Lonc and Elzobi \cite{LE99} that $\mathbf{m}\times \mathbf{n}$ is nice for any positive integers $m$ and $n$, it is very likely that the nice property holds at least for distributive lattices which are products of chains.
\begin{conjecture}\label{conjecture:prodchains}
For any positive integers $n_1,n_2,\ldots,n_r$, the lattice $\mathbf{n_1}\times \mathbf{n_2}\times \cdots\times \mathbf{n_r}$ is nice.
\end{conjecture}

Proposition \ref{prop-spos-nice} tells that a Schur positive graph must be nice. Together with Theorem \ref{main-thm-not-nice}, we know that not every distributive lattice is Schur positive. Motivated by Conjecture \ref{conjecture:StanleyBn} it would be very charming if a nice  distributive lattice is Schur positive. 
Unfortunately, this is still not true. In the next section, we will construct a family of nice distributive lattices which are not Schur positive.   

\section{Nice distributive lattices which are not Schur positive}  \label{section-4}

As we mentioned earlier, for any positive integers $m$ and $n$ Lonc and Elzobi \cite{LE99}  proved that $\mathbf{m}\times \mathbf{n}$ is nice. Somewhat surprisingly, the product $\mathbf{m}\times \mathbf{n}$ is not Schur postivie for many values of $m$ and $n$. 
In this paper, we denote the poset $\PP{m, n}$ as $\PP{(n+k), n}$, where $m=n+k$.
The main result of this section is as follows. 

\begin{theorem}\label{main-theorem-section4}
    For any positive integers $k\geq 5$ and $n\geq \frac{k+2}{2}$, 
    the distributive lattice $\PP{(n+k), n}$ 
    is not Schur positive.
\end{theorem}

Let us first outline the idea to prove Theorem \ref{main-theorem-section4}. It suffices to show that for some special partition $\rho$ the coefficient $[s_{\rho}]_{X_{\inc(\PP{(n+k), n})}}$ is negative. Here we shall take $\rho=(2n+k-1,2n+k-3,\ldots,k+3,k-3,2,2)$, a partition of $n(n+k)$ of length $n+2$. 
By Theorem \ref{theorem:WW} or Corollary \ref{coro:WW}, we need to enumerate elements $T\in\T_\rho$ such that $|\SCP_{\PP{(n+k), n},\cont(T)}|\neq 0$. Note that for any $T\in\T_\rho$ there are some restrictions on the content of $T$, say $\cont(T)=\la$.
On the one hand,  
$\la_1+\cdots+\la_i \geq \rho_1+\cdots+\rho_i$ for any integer $i$ since $\cont(T)\unrhd \sh(T)$; 
on the other hand, 
    in order to make $\SCP_{\PP{(n+k), n},\cont(T)}$ nonempty, $\la_1$ can be no larger than $\rho_1$ since a chain of $\PP{(n+k), n}$ contains at most $2n+k-1$ elements, thus $\la_1=\rho_1$; $\la_2$ can be no larger than $\rho_2$ since the next longest chain contains at most $2n+k-3$ elements, thus $\la_2=\rho_2$; continuing in this manner, we deduce that $\la_i=\rho_i$ for $i=1,\ldots,n-1$, and the first $n-1$ special rim hooks are exactly the first $n-1$ rows, as illustrated in \fref{bottomrimhooks}. Therefore, we have the following lemma. 

\begin{lemma}\label{lemma-content-speical}
For any $n,k\geq 1$, let $T$ be a special rim hook tabloid of shape $\rho=(2n+k-1,2n+k-3,\ldots,k+3,k-3,2,2)$ and content $\la$
such that $|\SCP_{\PP{(n+k), n},\la}|\neq 0$. 
Then we have $\la_i=\rho_i=2n+k-2i+1$ for $1\leq i\leq n-1$. 
\end{lemma}

    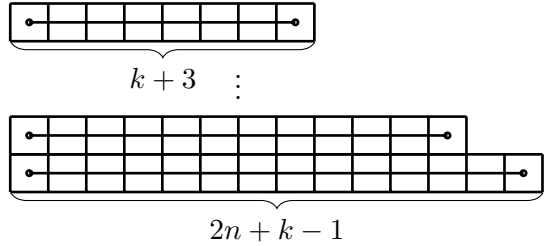
\begin{figure}[ht!]
        \centering
        \begin{tikzpicture}[scale=0.5]
            \ptnF{14,12}{0,0}
            \ptnF{8}{0,4}
            \specialribbonh{1,1}{13}
            \specialribbonh{1,2}{11}
            \specialribbonh{1,5}{7}
            \draw (6,3) node    {$\mathbf{\vdots}$};
            \draw[ decorate,decoration={brace,amplitude=6pt,mirror},xshift=0.4pt,yshift=-0.4pt] (0,0) -- (14,0) node[midway,yshift=-.5cm] {$2n+k-1$};
            \draw[ decorate,decoration={brace,amplitude=6pt,mirror},xshift=0.4pt,yshift=-0.4pt] (0,4) -- (8,4) node[midway,yshift=-.5cm] {$k+3$};
        \end{tikzpicture}
        \vspace*{-2mm}
        \caption{The first $n-1$ special rim hooks of any $T\in\T_{(2n+k-1,2n+k-3,\ldots,k+3,k-3,2,2)}$ satisfying $|\SCP_{\mathbf{(n+k)}\times \mathbf{n},\cont(T)}|\neq 0$.} \label{fig:bottomrimhooks}
    \end{figure}

Though for our purpose it would be enough to compute  
$|\SCP_{\PP{(n+k), n},\la}|$ with $\la_i=2n+k-2i+1$ for $1\leq i\leq n-1$ due to Lemma \ref{lemma-content-speical}, we shall give a more general result. For any $m\geq n\geq 1$, let $\la=(\la_1,\ldots,\la_\ell)$ be a partition of $m\times n$ with 
$\la_i=m+n-2i+1$ for $1\leq i\leq n-1$. One can check that $\la^{\langle n \rangle}$ is a partition of $m-n+1$, and 
write $\la^{\langle n \rangle}=\langle 1^{\al_1}, \ldots, (m-n+1)^{\al_{m-n+1}}\rangle$ with $\al_i$ being the number of occurrences of $i$ in $\la^{\langle n \rangle }$. The following result gives a formula to compute $|\SCP_{\mathbf{m}\times \mathbf{n},\la}|$ for any such partition $\la$.

\begin{proposition}\label{prop-main-SCP}
If $m\geq n\geq 1$, $\la=(\la_1,\ldots,\la_\ell)$ with $\la_i=m+n-2i+1$ for each $1\leq i\leq n-1$, and $\la^{\langle n \rangle}=\langle 1^{\al_1}, \ldots, (m-n+1)^{\al_{m-n+1}}\rangle$, then 
\begin{equation}\label{mnformula}
    |\SCP_{\mathbf{m}\times \mathbf{n},\la}| = (n-1)! \hspace{-3mm}
    \sum_{{(a_{i1},\,\ldots,\,a_{in})\vDash \al_i}\atop
    {1\leq i\leq m-n+1}} 
    \prod_{j=1}^{n}{\left(\sum_{k=1}^{m-n+1} k\cdot a_{kj}\right)!}  
    \prod_{k=1}^{m-n+1} \binom{\al_k}{a_{k1},\ldots,a_{kn}}\left(\frac{1}{k!}\right)^{\al_k}.
\end{equation}
\end{proposition}
\begin{proof} Let us first prove \eqref{mnformula} for small values of $n$.
When $n=1$  a little thought shows that 
\begin{align*}
    |\SCP_{\mathbf{m}\times \mathbf{n},\la}| = \binom{m}{\lambda_1,\,\lambda_2,\,\ldots,\,\lambda_{\ell}}=\frac{m!}{(1!)^{\alpha_1}\cdots (m!)^{\alpha_m}},    
\end{align*}
which coincides with \eqref{mnformula}.

We proceed to consider the case $n=2$. In this case, we need to enumerate the number of semi-ordered chain partitions of $\mathbf{m}\times \mathbf{2}$ of type $\la = (\la_1,\ldots,\la_\ell)$ with $\la_1 = m+1$. Note that for any $m\geq 2$ the maximal length that a chain of $\mathbf{m}\times \mathbf{2}$ can achieve is $m$ (containing $m+1$ elements), and  such a chain must connect $\hat{0}=(1,1)$ and $\hat{1}=(m,2)$, as illustrated in \fref{Pm2} for $m=8$ and $n=2$.  
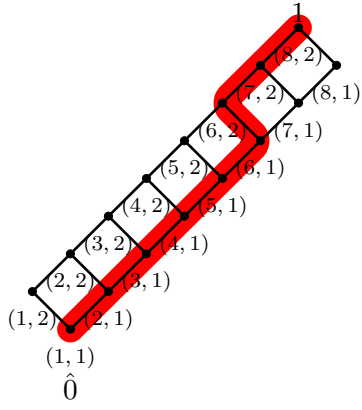
\begin{figure}[ht!]
    \centering
    \begin{tikzpicture}[scale=0.5]
        \colorchainl{0,0}{2,1,5}{red}
        \latticemn{8}{2}{0,0};
        \latticemnlabel{8}{2}{0,0};
        \node at (0,0.5) {$\hat{1}$};
        \node at (-6,-9.5) {$\hat{0}$};
    \end{tikzpicture}
    \caption{The poset $\mathbf{8}\times \mathbf{2}$ and a chain of size $9$.}\label{fig:Pm2}
\end{figure}

Looking at \fref{Pm2}, we see that $\mathbf{m}\times \mathbf{2}$ has $m$ chains of size $m+1$, all connecting $\hat{0}$ and $\hat{1}$. Moreover, each chain $C$ with $m+1$ elements is determined by an $r\in [m]=\{1,\ldots,m\}$ such that $\{(r,1),(r,2)\}\subseteq C$. 
Once the chain $C$ is removed, the remaining elements of $\mathbf{m}\times \mathbf{2}$ will form $2$ chains, say $R_1=\{(1,2),(2,2),\ldots,(r-1,2)\}$ and $R_2=\{(r+1,1),\ldots,(m,1)\}$ (each may be empty).
We further need to partition the set $R_1\cup R_2$ into chains whose sizes form a partition $\la^{\langle 2 \rangle}=(\la_2,\ldots,\la_\ell)$. 
Since any two elements $a\in R_1$ and $b\in R_2$ are incomparable, each chain in a chain partition of $R_1\cup R_2$  
lies entirely inside $R_1$ or entirely inside $R_2$. Recall that $\la^{\langle 2 \rangle}=\langle 1^{\al_1}, \ldots, (m-1)^{\al_{m-1}}\rangle$. 
If we fix a weak composition $(a_{k1},a_{k2})\vDash \alpha_k$ for each $1\leq k\leq m-1$, then the number of semi-ordered chain partitions of $R_1\cup R_2$ of type $\la^{\langle 2 \rangle}$, such that there are $a_{k1}$ chains of size $k$ in $R_1$ and $a_{k2}$ chains of size $k$ in $R_2$, is equal to
\begin{equation*}
    \frac{(r-1 )!}{\prod_{k=1}^{m-1}(k!)^{a_{k1}}\cdot a_{k1}!} \cdot \frac{(m-r) !}{\prod_{k=1}^{m-1}(k!)^{a_{k2}}\cdot a_{k2}!} \cdot \prod_{k=1}^{m-1} \al_k!
    =
    (r-1)!(m-r)! \prod_{k=1}^{m-1} \binom{\al_k}{a_{k1},a_{k2}}\left(\frac{1}{k!}\right)^{\al_k}.
\end{equation*}
Since $|R_1|=r-1=\sum_{k=1}^{m-1} k\cdot a_{k1}$ and $|R_2|=m-r=\sum_{k=1}^{m-1} k\cdot a_{k2}$,
we have
\begin{align*}
   |\SCP_{\mathbf{m}\times \mathbf{2},\la}| &= \sum_{(a_{i1},a_{i2})\vDash\al_i\atop 1\leq i\leq m-1} 
   \left(\sum_{k=1}^{m-1} k\cdot a_{k1}\right)!\left(\sum_{k=1}^{m-1} k\cdot a_{k2}\right)! \prod_{k=1}^{m-1} \binom{\al_k}{a_{k1},a_{k2}}\left(\frac{1}{k!}\right)^{\al_k},
\end{align*}
as desired. 

Now we move to the proof for general $n$, and the main idea is the same as the case $n=2$. 
To determine a semi-ordered chain partition of $\mathbf{m}\times \mathbf{n}$ with the desired type 
$\la$, we first take $n-1$ disjoint chains whose sizes are $\la_1,\ldots,\la_{n-1}$ respectively. 
As will be shown later, once these $n-1$ chains are taken out from $\mathbf{m}\times \mathbf{n}$, the remaining part can be written as a disjoint union of at most $n$ chains, each of which is a saturated chain of the subposet consisting of 
$\{(1,k),(2,k),\ldots,(m,k)\}$ for some $k$.  

In order to choose the first $n-1$ chains, let us first enumerate 
the number of elements for each rank in the poset $\mathbf{m}\times \mathbf{n}$.
Note that, for each $1\leq k\leq n-1$, there are $k$ elements of rank $k-1$ and rank $m+n-k-1$ respectively, and, for each $n-1\leq i\leq m-1$, there are $n$ elements of rank $i$. 
For example, the poset $\mathbf{16}\times \mathbf{4}$ in \fref{Pmn} has $1$ element of rank $0$ and rank $18$, $2$ elements of rank $1$ and rank $17$, $3$ elements of rank $2$ and rank $16$, and $4$ elements of rank $i$ for $i=3,\ldots,15$. 

Keep in mind that we want to enumerate semi-ordered chain partitions of $\mathbf{m}\times \mathbf{n}$ of type $\la = (\la_1,\ldots,\la_\ell)$ with $\la_1 = m+n-1,\ \la_2 = m+n-3,\ \ldots,\ \la_{n-1} = m-n+3$.
We will choose the first $n-1$ chains according to their sizes one by one. The first chain is of size $\la_1=m+n-1$, which will use up all elements of ranks $0$ and $m+n-2$, and also use one element of rank $i$ for each $1\leq i\leq m+n-3$. Once such a chain is taken from  $\mathbf{m}\times \mathbf{n}$, the second chain of size $\la_2 = m+n-3$, if it exists, will use up the remaining elements of ranks $1$ and $m+n-3$, and also use one element of rank $i$ for each $2\leq i\leq m+n-4$. Continuing in this manner, 
the $k$-th chain of size $\la_{k} = m+n-2k+1$, if it exists, will use up the remaining elements of ranks $k-1$ and $m+n-k-1$, and also use one element of rank $i$ for each $k\leq i\leq m+n-k-2$. Finally, 
the $(n-1)$-th chain of size $\la_{n-1} = m-n+3$, if it exists, will use up the remaining elements of ranks $n-2$ and $m$, and also use one element of rank $i$ for each $n-1\leq i\leq m-1$. For an illustration of the above procedure, see \fref{Pmn}. 

Note that such a choice of the above $n-1$ chains will use up all elements of ranks $0,1,\ldots,n-2$ and $m,m+1,\ldots,m+n-2$, leaving only one element for each rank in $\{n-1,\ldots,m-1\}$. 
Due to their saturated property that each element $(i,j)$ connects either $(i+1,j)$ or $(i,j+1)$, these $n-1$ chains between rank $n-2$ and rank $m$ form a non-crossing path family from $\{(1,n-1),(2,n-2),\ldots,(n-1,1)\}$ to $\{(m-n+2,n),(m-n+3,n-1),\ldots,(m,2)\}$ in the Hasse diagram of $\mathbf{m}\times \mathbf{n}$.
Since there are exactly $n-1$ elements of rank $n-2$, 
one can associate a permutation $\sigma\in \mathfrak{S}_{n-1}$ to these first $n-1$ chains by letting the $i$-th chain from left to right (which contains $(i,n-i)$ and $(m-n+1+i,n+1-i)$ simultaneously) to be of size $\la_{\sigma_i}$.  
Now the set of the remaining elements is of the form $\{(1,n),(2,n),\ldots,(r_1-1,n),(r_1+1,n-1),\ldots,(r_2-1,n-1),\ldots,(r_{n-1}+1,1),\ldots,(m,1)\}$, determined by some $(n-1)$ elements $r_1<r_2<\cdots<r_{n-1}$ of $[m]$. Set $r_0=0$ and $r_n=m+1$ and for $1\leq i \leq n$ let 
\begin{align}\label{eq-R}
R_i=\{(r_{i-1}+1,n+1-i),\ldots,(r_{i}-1,n+1-i)\}.
\end{align}
For example, the semi-ordered chain partition in \fref{Pmn} has $R_1=\{(1,4),(2,4),(3,4)\},R_2=\{(5,3),(6,3),(7,3),(8,3),(9,3)\},R_3=\{(11,2),(12,2),(13,2)\}$ and $R_4=\{(15,1),(16,1)\}$.
We see that each $R_i$ is a saturated chain, and their union is just the set of the remaining elements after taking the first $n-1$ chains from $\mathbf{m}\times \mathbf{n}$. 
To summarize, a collection $\mathcal{C}$ of $n-1$ chains of respective sizes $\la_1,\ldots,\la_{n-1}$ will determine a permutation $\sigma_{\mathcal{C}}\in \mathfrak{S}_{n-1}$ and an $(n-1)$-element set $\Upsilon_{\mathcal{C}}=\{r_1,\ldots,r_{n-1}\}\subseteq [m]$. In fact, the converse is also true. 

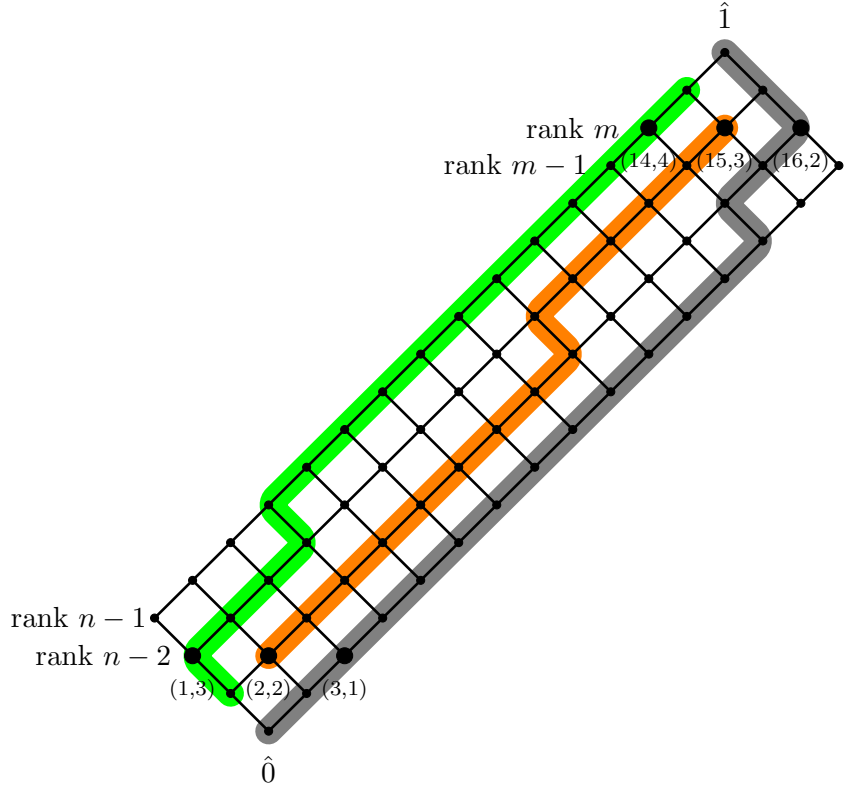
\begin{figure}[ht!]
    \centering
    \begin{tikzpicture}[scale=0.5]
        \colorchain{0,0}{2,2,1,13}{gray}
        \colorchainl{-1,-1}{11,1,3,1}{green}
        \colorchainl{0,-2}{5,1,8}{orange}
        \latticemn{16}{4}{0,0};
        \latticemnlabelrank{16}{4}{0,0}{2};
        \latticemnlabelrank{16}{4}{0,0}{16};
        \node at (0,1) {$\hat{1}$};
        \node at (-12,-19) {$\hat{0}$};
        \node at (-4,-2) {rank $m$};
        \node at (-5.5,-3) {rank $m-1$};
        \node at (-17,-15) {rank $n-1$};
        \node at (-16.35,-16) {rank $n-2$};
    \end{tikzpicture}
    \caption{The poset $\mathbf{16}\times \mathbf{4}$.}\label{fig:Pmn}
\end{figure}

{\bf (Claim I) } Any
permutation $\sigma\in \mathfrak{S}_{n-1}$ and any $(n-1)$-element set $\Upsilon=\{r_1,\ldots,r_{n-1}\}\subseteq [m]$ also determine a collection $\mathcal{C}$ of $n-1$ chains of sizes $\la_1,\ldots,\la_{n-1}$ such that $\sigma=\sigma_{\mathcal{C}}$ and $\Upsilon=\Upsilon_{\mathcal{C}}$. 

Now let us give a proof of this claim. Our aim is to construct a collection $\mathcal{C}$ of $n-1$ chains $C_1, C_2,\ldots,C_{n-1}$
such that $C_i$ is of size $\la_i=m+n-2i+1$. 
Suppose that $r_1<r_2<\cdots<r_{n-1}$, let us define $R_1,\ldots,R_{n-1}$ by using \eqref{eq-R}.
For each $1\leq i\leq n-1$ let $C_{\sigma_i}$ 
contain the elements $\{(i,n-i),\ldots,(r_i,n-i),(r_i,n-i+1),\ldots,(m-n+1+i,n-i+1)\}$.
Thus the $k$-th element of rank $n-2$ from left to right is assigned to the chain $C_{\sigma_k}$
for each $1\leq k\leq n-1$. Now the chain $C_{n-1}$ has been constructed. 
To complete the proof, it suffices to show how to use $\sigma$ to assign elements of ranks $0,\ldots,n-3$ and  $m+1,\ldots,m+n-2$ to $C_1, C_2,\ldots,C_{n-2}$. 
For each $1\leq i\leq n-2$ let $\sigma^{(i)}$ denote the permutation obtained from $\sigma$ by deleting all numbers greater than $i$.
Now, for each $1\leq k\leq i$, assign the $k$-th element of rank $i-1$ and the $k$-th element of rank $m+n-i-1$ from left to right to the chain $C_{\sigma^{(i)}_k}$. Now we construct all the chains $C_1, C_2,\ldots,C_{n-1}$. 
It is routine to check that each $C_{i}$ is from an element of rank $i-1$ to an element of rank $m+n-i-1$,
which has the desired size $\la_i$. Moreover, we have $\sigma=\sigma_{\mathcal{C}}$ and $\Upsilon=\Upsilon_{\mathcal{C}}$. This completes the proof of the claim. 

\fref{Pmnfactorial} gives an example of the assignment of elements of ranks $0,1,\ldots,6$ when $n=8$ and $\sg=3165274\in\Sn{7}$. The elements of rank $6$ will be assigned to chains $C_3,C_1,C_6,C_5,C_2,C_7,C_4$ from left to right; the elements of rank $5$ will be assigned to chains $C_3,C_1,C_6,C_5,C_2,C_4$ from left to right; the elements of rank $4$ will be assigned to chains $C_3,C_1,C_5,C_2,C_4$ from left to right. Continuing in this manner, all the $28$ elements of ranks $0,1,\ldots,6$ are assigned to the $7$ chains as shown in \fref{Pmnfactorial}.

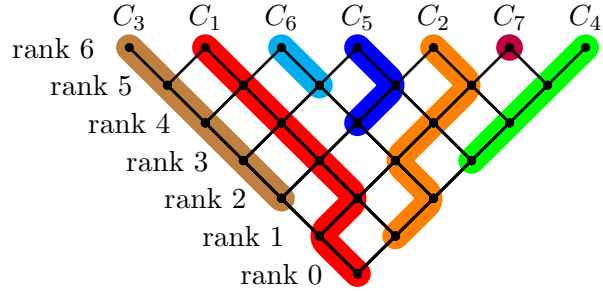
\begin{figure}[ht!]
    \centering
    \begin{tikzpicture}[scale=0.5]
        
        \fill[opacity=0.5,purple] (4,6) circle (10pt);		
        \colorchain{-4,6}{4,1,1}{red}
        \colorchain{2,6}{1,2,1,1}{orange}
        \colorchain{-6,6}{4}{brown}
        \colorchainl{6,6}{3}{green}
        \colorchain{0,6}{1,1}{blue}
        \colorchain{-2,6}{1}{cyan}
        \latticemnud{1}{7}{0,0};
        \latticemnud{2}{6}{0,0};
        \latticemnud{3}{5}{0,0};
        \latticemnud{4}{4}{0,0};
        \latticemnud{5}{3}{0,0};
        \latticemnud{6}{2}{0,0};
        \latticemnud{7}{1}{0,0};
        \node at (-6,6.8) { $C_3$};
        \node at (-4,6.8) { $C_1$};
        \node at (-2,6.8) { $C_6$};
        \node at (0,6.8) { $C_5$};
        \node at (2,6.8) { $C_2$};
        \node at (4,6.8) { $C_7$};
        \node at (6,6.8) { $C_4$};
        \node at (-8,6) {rank $6$};
        \node at (-7,5) {rank $5$};
        \node at (-6,4) {rank $4$};
        \node at (-5,3) {rank $3$};
        \node at (-4,2) {rank $2$};
        \node at (-3,1) {rank $1$};
        \node at (-2,0) {rank $0$};
    \end{tikzpicture}
    \caption{An assignment of elements of first 7 ranks of $\mathbf{m}\times\mathbf{8}$.}\label{fig:Pmnfactorial}
\end{figure}


Given an $(n-1)$-subset $\Upsilon=\{r_1,\ldots,r_{n-1}\}$ of $[m]$, let
$R_{\Upsilon}$ denote the subposet induced by $\cup_{i=1}^{n}R_i$, where $R_i$ is defined by \eqref{eq-R}.
To complete the proof of the proposition, we further need to enumerate the semi-ordered chain partitions of $R_{\Upsilon}$ of type $\la^{\langle n \rangle}=(\la_n,\ldots,\la_\ell)$.
Notice that each $R_i$ form a chain (may be empty), and any two elements in different $R_i$'s are incomparable. Hence all elements of every chain in the chain partition belong to the same $R_i$.
Thus,  the number of semi-ordered chain partitions of $R_{\Upsilon}$ of type $\la^{\langle n \rangle}=\langle 1^{\al_1},\,\ldots,\,(m-n+1)^{\al_{m-n+1}}\rangle$, such that $a_{kj}$ chains of size $k$ are picked from $R_j$, is equal to
\begin{equation*}
    \prod_{j=1}^{n}\frac{(\sum_{k=1}^{m-n+1} k\cdot a_{kj})!}{\prod_{k=1}^{m-n+1}(k!)^{a_{kj}}\cdot a_{kj} !} 
    \cdot \prod_{k=1}^{m-n+1} \al_k!
    =
    \prod_{j=1}^{n}{\left(\sum_{k=1}^{m-n+1} k\cdot a_{kj}\right)!}  
    \prod_{k=1}^{m-n+1} \binom{\al_k}{a_{k1},\ldots,a_{kn}}\left(\frac{1}{k!}\right)^{\al_k}.
\end{equation*}
Then
\begin{align}
|\SCP_{R_{\Upsilon},\la^{\langle n \rangle}}|= \sum_{(a_{i1},\ldots,a_{in})} 
    \prod_{j=1}^{n}{\left(\sum_{k=1}^{m-n+1} k\cdot a_{kj}\right)!}  
    \prod_{k=1}^{m-n+1} \binom{\al_k}{a_{k1},\ldots,a_{kn}}\left(\frac{1}{k!}\right)^{\al_k},\label{eq-temp-mnformula}
\end{align}
summing over all weak compositions $(a_{i1},\ldots,a_{in}) \vDash \al_i$ such that $\sum_{k=1}^{m-n+1} k\cdot a_{kj}=|R_j|$ for each $1\leq i\leq m-n+1$ and $1\leq j\leq n$. 
From Claim I it follows that
\begin{align}\label{eq-temp1-mnformula}
|\SCP_{\mathbf{m}\times \mathbf{n},\la}|&=\sum_{\sigma\in \mathfrak{S}_{n-1}}\sum_{\Upsilon\subseteq [m]}
|\SCP_{R_{\Upsilon},\la^{\langle n \rangle}}|=(n-1)!\sum_{\Upsilon\subseteq [m]}
|\SCP_{R_{\Upsilon},\la^{\langle n \rangle}}|,
\end{align}
where the second equality holds since $|\SCP_{R_{\Upsilon},\la^{\langle n \rangle}}|$ is independent of the permutation $\sigma$.
Combining \eqref{eq-temp-mnformula}  and \eqref{eq-temp1-mnformula} leads to \eqref{mnformula}, as desired.  This completes the proof. 
\end{proof}

\begin{remark}
In view of our formula \eqref{mnformula}, the number $|\SCP_{\mathbf{m}\times \mathbf{n},\la}|$ can also be computed by the following two formulas: 
\begin{equation*}
    |\SCP_{\mathbf{m}\times \mathbf{n},\la}| = (n-1)! \left(\prod_{k=1}^{m-n+1} (\sum_{i=1}^{n}\frac{x_i^k}{k!})^{\al_k} \right)
    \left(\sum_{(\ga_1,\ldots,\ga_n)\vDash m-n+1} \prod_{j=1}^{n}\frac{\ga_j!}{x_j^{\ga_j}}\right)\Bigg|_\textnormal{constant term},
\end{equation*}
and
\begin{equation*}
    |\SCP_{\mathbf{m}\times \mathbf{n},\la}| = (n-1)!  \left(\sum_{(\ga_1,\ldots,\ga_n)\vDash m-n+1} 
    {\partial {\bf x}^\ga}\right)
    \left(\prod_{k=1}^{m-n+1} \bigg(\sum_{i=1}^{n}\frac{x_i^k}{k!}\bigg)^{\al_k} \right),
\end{equation*}
where $\partial{\bf x}^\ga=\partial x_1^{\ga_1}\cdots \partial x_n^{\ga_n}$.
\end{remark}

Now, we are ready to prove Theorem \ref{main-theorem-section4}.

\textit{Proof of Theorem \ref{main-theorem-section4}.}
    As discussed in the paragraph immediately after Theorem \ref{main-theorem-section4}, 
    it suffices to show that the coefficient $[s_{\rho}]_{X_{\inc(\PP{(n+k), n})}}$ is negative for $\rho=(2n+k-1,2n+k-3,\ldots,k+3,k-3,2,2)$. 
    In order to use (\ref{formula}) to compute this coefficient, we need to enumerate 
    elements $T\in\T_\rho$ such that $|\SCP_{\PP{(n+k), n},\cont(T)}|\neq 0$. 
    By Lemma \ref{lemma-content-speical}, if $T\in\T_\rho$ has content $\la$, then 
    $\la_i=\rho_i$ for $1\leq i\leq n-1$ and the first $n-1$ special rim hooks are exactly the first $n-1$ rows, as shown in \fref{bottomrimhooks}.
    We continue to consider the special rim hooks in the remaining top $3$ rows of $T\in\T_\rho$. 
    Due to its special shape, $T$ has only $6$ possibilities, labeled as $T_1,\,T_2,\,\ldots,\,T_6$ in \fref{toprimhooks}, where for each $T_i$ only special rim hooks in the top $3$ rows are shown. 
    Thus $|\T_{\rho}|=6,\, \wt(T_1)=3,\, \wt(T_2)=\wt(T_3)=2,\, \wt(T_4)=\wt(T_5)=1$ and
    $\wt(T_6)=0$. We further need to compute $|\SCP_{\PP{(n+k),n},\cont(T_i)}|$ for each $1\leq i\leq 6$, 
    which can be done by using Proposition \ref{prop-main-SCP}. 
    For convenience of notation, set $\delta_{n,k}=(2n+k-1,2n+k-3,\ldots,k+3)$. 
    
    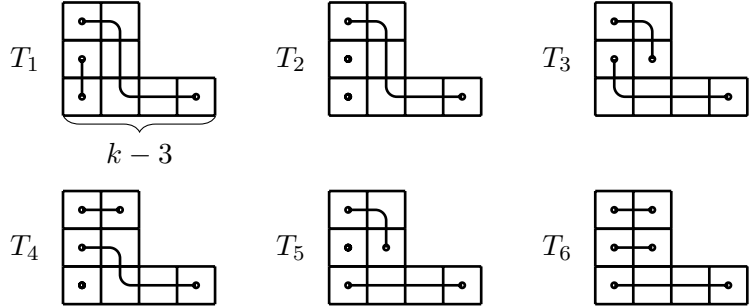
\begin{figure}[ht!]
        \centering
        \begin{tikzpicture}[scale=0.5]
            \ptnF{4,2,2}{0,0}
            \specialribbonh{1,3}{1,2,2}
            \specialribbonv{1,2}{1}
            \draw (-1,1.5) node {$T_1$};
            \draw[decorate,decoration={brace,amplitude=6pt,mirror},xshift=0.4pt,yshift=-0.4pt] (0,0) -- (4,0) node[midway,yshift=-.5cm] {$k-3$};
            
            \ptnF{4,2,2}{7,0}
            \specialribbonh{8,3}{1,2,2}
            \specialribbonv{8,2}{0}
            \specialribbonv{8,1}{0}
            \draw (6,1.5) node {$T_2$};
            \ptnF{4,2,2}{14,0}
            \specialribbonh{15,3}{1,1}
            \specialribbonv{15,2}{1,3}
            \draw (13,1.5) node {$T_3$};
            \ptnF{4,2,2}{0,-5}
            \specialribbonh{1,-2}{1}
            \specialribbonh{1,-3}{1,1,2}
            \specialribbonv{1,-4}{0}
            \draw (-1,-3.5) node {$T_4$};
            \ptnF{4,2,2}{7,-5}
            \specialribbonh{8,-2}{1,1}
            \specialribbonh{8,-3}{0}
            \specialribbonh{8,-4}{3}
            \draw (6,-3.5) node {$T_5$};
            \ptnF{4,2,2}{14,-5}
            \specialribbonh{15,-2}{1}
            \specialribbonh{15,-3}{1}
            \specialribbonh{15,-4}{3}
            \draw (13,-3.5) node {$T_6$};
        \end{tikzpicture}
        \caption{The special rim hooks in top $3$ rows of $T_i\in\T_{(2n+k-1,2n+k-3,\ldots,k+3,k-3,2,2)}$.}\label{fig:toprimhooks}
    \end{figure}

    Let us first compute $|\SCP_{\PP{(n+k),n},\cont(T_1)}|$. Note that $\cont(T_1)=(\delta_{n,k},k-1,2)$, and hence $\cont(T_1)^{\langle n\rangle}=\langle 2^1, (k-1)^1\rangle$ since $k\geq 5$. Based on \eqref{mnformula}, we obtain 
        \begin{align*}
        |\SCP_{\PP{(n+k),n},\cont(T_1)}|=&(n-1)!\times \sum_{\substack{(a_{1},\ldots,a_{n})\vDash 1  \\ (b_{1},\ldots,b_{n})\vDash 1}}
        \prod_{j=1}^{n}{(2\cdot a_{j} + (k-1)\cdot b_{j})!} \times 
        \frac{1}{(k-1)!\cdot 2!}\\[8pt]
        =& (n-1)!\times \frac{n\cdot (k+1)!+n(n-1)\cdot (k-1)!\cdot 2!}{(k-1)!\cdot 2!}\\[6pt]
        =& n!\left(n+\frac{k^2+k-2}{2}\right).
    \end{align*}

    For $T_2$ we have $\cont(T_2)=(\delta_{n,k},k-1,1,1)$ and hence $\cont(T_2)^{\langle n\rangle}=\langle 1^2,(k-1)^1 \rangle$. From \eqref{mnformula} it follows that
    \begin{align*}
        |\SCP&_{\PP{(n+k),n},\cont(T_2)}|\\[6pt]
        =&(n-1)!\times \sum_{\substack{(a_{1},\ldots,a_{n})\vDash 2  \\ (b_{1},\ldots,b_{n})\vDash 1}}
        \prod_{j=1}^{n}{(1\cdot a_{j} + (k-1)\cdot b_{j})!} \cdot 
        \binom{2}{a_1,\ldots,a_n}\frac{1}{(k-1)!\cdot (1!)^2}\\[8pt]
        =&(n-1)!\times \frac{n(n-1)(n-2)\cdot (k-1)!+2n(n-1)\cdot k!+2n(n-1)\cdot (k-1)!+n\cdot (k+1)!}{(k-1)!\cdot (1!)^2}\\[6pt]
        =&n!(n^2+(2k-1)n+(k^2-k)).
    \end{align*}
    
    For $T_3$ one can check that $\cont(T_3)=(\delta_{n,k},k-2,3)$.  
    When $k=5$, we have $\cont(T_3)^{\langle n\rangle}=\langle 3^2\rangle$. By using \eqref{mnformula}, we find that
    \begin{align*}
    |\SCP_{\PP{(n+5),n},\cont(T_3)}|
    =&(n-1)!\times\sum_{(a_{1},\ldots,a_{n})\vDash 2}
        \prod_{j=1}^{n}{(3\cdot a_{j})!} \cdot \binom{2}{a_1,\ldots,a_n}\cdot
        \frac{1}{(3!)^2}\\[8pt]
        =&(n-1)!\times \frac{n \cdot 6!+n(n-1)\cdot 3!\cdot 3!}{(3!)^2}\\[6pt]
        =&n!(n+19).
    \end{align*}
    When $k\geq 6$, we have $\cont(T_3)^{\langle n\rangle}=\langle 3^1,(k-2)^1\rangle$. Similarly, by \eqref{mnformula}
    we have
    \begin{align*}
    |\SCP_{\PP{(n+k),n},\cont(T_3)}|
        =&(n-1)!\times \sum_{\substack{(a_{1},\ldots,a_{n})\vDash 1  \\ (b_{1},\ldots,b_{n})\vDash 1}}
        \prod_{j=1}^{n}{(3\cdot a_{j} + (k-2)\cdot b_{j})!} \cdot 
        \frac{1}{(k-2)!\cdot 3!}\\[8pt]
        =&(n-1)!\times \frac{n\cdot (k+1)!+n(n-1)\cdot (k-2)!\cdot 3!}{(k-2)!\cdot 3!}\\[6pt]
        =&n!\left(n+\frac{k^3-k-6}{6}\right).
    \end{align*}
    
    For $T_4$ we notice that $\cont(T_4)=(\delta_{n,k},k-2,2,1)$ and hence $\cont(T_4)^{\langle n\rangle}= \langle 1^1, 2^1, (k-2)^1 \rangle$
    for any $k\geq 5$. By using \eqref{mnformula} again, we obtain
    \begin{small}   
    \begin{align*}
        |\SCP&_{\PP{(n+k),n},\cont(T_4)}|\\
        =&(n-1)!\times \sum_{\substack{(a_{1},\ldots,a_{n})\vDash 1  \\ (b_{1},\ldots,b_{n})\vDash 1 \\ (c_{1},\ldots,c_{n})\vDash 1}}
        \prod_{j=1}^{n}{(1\cdot a_{j} + 2\cdot b_{j} + (k-2)\cdot c_{j})!} \cdot 
        \frac{1}{(k-2)!\cdot 2!\cdot 1!}\\[8pt]
        =&(n-1)!\times \frac{n(n-1)(n-2)\cdot (k-2)!\cdot 2!+n(n-1)(k!+(k-1)!\cdot 2!+(k-2)!\cdot 3!)+n\cdot (k+1)!}{2\cdot (k-2)!}\\[6pt]
        =&n!\left(n^2+\frac{k^2+k-2}{2}n+\frac{k^3-k^2-2k}{2}\right).
    \end{align*}
    \end{small}\
    Note that the tabloid $T_5$ has content $\cont(T_5)=(\delta_{n,k},k-3,3,1)$ for $k\geq 5$. 
    When $k=6$ we have $\cont(T_5)^{\langle n\rangle}=\langle 1^1, 3^2 \rangle $. Applying \eqref{mnformula} gives
    \begin{align*}
        |\SCP&_{\PP{(n+6),n},\cont(T_5)}|\\
        =&(n-1)!\times \sum_{\substack{(a_{1},\ldots,a_{n})\vDash 1  \\ (b_{1},\ldots,b_{n})\vDash 2 }}
        \prod_{j=1}^{n}{(1\cdot a_{j} + 3\cdot b_{j})!} \cdot 
        \binom{2}{b_{1},\ldots,b_{n}}\frac{1}{(3!)^2\cdot 1!}\\[8pt]
        =&(n-1)!\times \frac{n(n-1)(n-2)\cdot (3!)^2\cdot 1! +n(n-1)\cdot 6!\cdot 1!+2n(n-1)\cdot 3!\cdot 4!+ n\cdot 7!}{(3!)^2 1!}\\[6pt]
        =&n!(n^2+25n+114).
    \end{align*}
    It is also clear that $\cont(T_5)^{\langle n\rangle}= \langle 1^1, (k-3)^1, 3^1 \rangle$  for $k=5$, and $\cont(T_5)^{\langle n\rangle}= \langle 1^1, 3^1, (k-3)^1 \rangle$ for $k\geq 7$. Thus, when $k\neq 6$, the formula \eqref{mnformula} gives
    \begin{small}
    \begin{align*}
        |\SCP&_{\PP{(n+k),n},\cont(T_5)}|\\
        =&(n-1)!\times \sum_{\substack{(a_{1},\ldots,a_{n})\vDash 1  \\ (b_{1},\ldots,b_{n})\vDash 1 \\ (c_{1},\ldots,c_{n})\vDash 1}}
        \prod_{j=1}^{n}{(1\cdot a_{j} + 3\cdot b_{j} + (k-3)\cdot c_{j})!} \cdot 
        \frac{1}{(k-3)!\cdot 3!\cdot 1!}\\[8pt]
        =&(n-1)!\times \frac{n(n-1)(n-2)\cdot (k-3)!\cdot 3!+n(n-1)(k!+(k-2)!\cdot 3!+(k-3)!\cdot 4!)+n(k+1)!}{(k-3)!\cdot 3!\cdot 1!}\\[6pt]
        =&n!\left(n^2+\frac{k^3-3k^2+8k-6}{6}n+\frac{k^4-3k^3+2k^2-6k}{6}\right).
    \end{align*}
    \end{small}
        
    For the tabloid $T_6$ we have $\cont(T_6)=(\delta_{n,k},k-3,2,2)$. 
    If $k\geq 6$, then $\cont(T_6)^{\langle n\rangle}= \langle 2^2, (k-3)^1 \rangle$. By using \eqref{mnformula} we get
    \begin{small}
    \begin{align*}
        |\SCP&_{\PP{(n+k),n},\cont(T_6)}|\\[5pt]
        =&(n-1)!\times \sum_{\substack{(a_{1},\ldots,a_{n})\vDash 2  \\ (b_{1},\ldots,b_{n})\vDash 1 }}
        \prod_{j=1}^{n}{(2\cdot a_{j} + (k-3)\cdot b_{j})!} \cdot 
        \binom{2}{a_{1},\ldots,a_{n}}\frac{1}{(k-3)!\cdot (2!)^2}\\[8pt]
        =&(n-1)!\times \frac{n(n-1)(n-2)(k-3)!(2!)^2+2n(n-1)(k-1)!2!+n(n-1)(k-3)!4!+n(k+1)!}{(k-3)!(2!)^2}\\[6pt]
        =&n!\left(n^2+(k^2-3k+5)n+\frac{k^4-2k^3-5k^2+14k-24}{4}\right).
    \end{align*}
    \end{small}
    When $k=5$, we have $\cont(T_6)^{\langle n\rangle}=\langle 2^3 \rangle$. Again using \eqref{mnformula} we obtain
    \begin{align*}
    |\SCP_{\PP{(n+5),n},\cont(T_6)}|=&(n-1)!\times \sum_{(a_{1},\ldots,a_{n})\vDash 3}
        \prod_{j=1}^{n}{(2\cdot a_{j})!} \cdot 
        \binom{3}{a_{1},\ldots,a_{n}}\frac{1}{(2!)^3}\\[8pt]
        =&(n-1)!\times \frac{n(n-1)(n-2)(2!)^3+3n(n-1)4!2!+6!n}{(2!)^3}\\
        =&n!(n^2+15n+74).
    \end{align*}
        
    Based on \eqref{formula} and the above enumerative results on $|\SCP_{\PP{(n+k),n},\cont(T_i)}|$, one can check that 
    for $\rho=(2n+k-1,2n+k-3,\ldots,k+3,k-3,2,2)$, 
    \begin{align*}
        [s_\rho]_{X_{\inc(\PP{(n+5),n})}}&=n!(-4n+9),\\
        [s_\rho]_{X_{\inc(\PP{(n+6),n})}}&=n!(-11n+32),
    \end{align*}
    which are both negative when $n\geq 3$. 
    Similarly, for $k\geq 7$, one can check that 
    \begin{align*}
        [s_\rho]_{X_{\inc(\PP{(n+k),n})}} &=\frac{n!(k-4)}{12} \left((-2 k^2 + 4 k - 18) n+(k^3 - 7 k + 18)\right).
    \end{align*}
    Now it suffices to show that for any $k\geq 7$ and $n\geq \frac{k+2}{2}$ the number $(-2 k^2 + 4 k - 18) n+(k^3 - 7 k + 18)$ 
    is negative. This is obvious since for $k\geq 7$ there holds
$2k^2-4k+18>0$ and 
\begin{equation*}
        \frac{k^3-7k+18}{2k^2-4k+18}<\frac{k+2}{2}. 
    \end{equation*}
This completes the proof.
\qed

We proceed to give a more general result than Theorem \ref{main-theorem-section4}. 
Given two posets $P$ and $Q$, let $P\oplus Q$ denote their direct sum. Based on the work of Lonc and Elzobi \cite{LE99}, we 
can obtain the following result. 

\begin{theorem}\label{thm-nice-dl}
For any $p,q\geq 0$ and $m,n\geq 1$, the poset
$\mathbf{p}\oplus (\mathbf{m}\times \mathbf{n})\oplus \mathbf{q}$ is a nice distributive lattice.
\end{theorem}

\begin{proof}
It is routine to check that $\mathbf{p}\oplus (\mathbf{m}\times \mathbf{n})\oplus \mathbf{q}$ is a distributive lattice.
It remains to show that this lattice is nice. 
Define a partition $\lambda=(\lambda_1,\ldots,\lambda_n)$ by letting $\lambda_i=m+n-2i+1$ for $1\le i\le n$.
Theorem 1 of Lonc and Elzobi \cite{LE99} tells that there exists a chain partition of 
$\mathbf{m}\times \mathbf{n}$ of type $\delta$ if and only if $\delta\unlhd\lambda$.
Set $\tilde{\lambda}=(\lambda_1+p+q,\lambda_2,\ldots,\lambda_n)$. 
A little thought shows that if there exists a chain partition of $\mathbf{p}\oplus (\mathbf{m}\times \mathbf{n})\oplus \mathbf{q}$
of type $\mu$, then  $\mu\unlhd \tilde{\lambda}$.
Thus, to prove that $\mathbf{p}\oplus (\mathbf{m}\times \mathbf{n})\oplus \mathbf{q}$ is nice, it suffices to show the converse of the above statement is also true.  
Suppose that $\mu=(\mu_1,\ldots,\mu_\ell)\unlhd \tilde{\lambda}$ and hence $\ell\geq n$. We define $t_1,\ldots,t_{\ell}$ recursively as follows. 
If $\mu_1\leq \lambda_1$, then set $t_1=0$, and otherwise, set $t_1=\mu_1-\lambda_1$.
Assuming that $t_1,\ldots,t_{j-1}$ have been defined, we set
\begin{align}\label{eq-tj}
t_j=\max\left\{0, \sum_{i=1}^j \mu_i-\sum_{i=1}^j\lambda_i-\sum_{i=1}^{j-1} t_i\right\}.
\end{align}
It is routine to check that
\begin{align}\label{eq-tj-1}
\sum_{i=1}^j \mu_i-\sum_{i=1}^{j} t_i\leq \sum_{i=1}^j\lambda_i, \quad \mbox{for any $1\leq j\leq \ell$.}
\end{align}
It follows that $t_j\geq 0$ and  
\[\sum_{i=1}^\ell t_i\ge \sum_{i=1}^\ell \mu_i-\sum_{i=1}^\ell \lambda_i=p+q.\]
We claim that $\sum_{i=1}^\ell t_i=p+q$. Suppose that $\sum_{i=1}^\ell t_i>p+q$. 
Let $j$ be the minimum integer such that $\sum_{i=1}^j t_i>p+q$, and hence $t_j\ge 1$. 
Thus by \eqref{eq-tj} we get
\[t_j=\sum_{i=1}^j \mu_i-\sum_{i=1}^j\lambda_i-\sum_{i=1}^{j-1} t_i.\]
On the other hand, from $\mu\unlhd\tilde{\lambda}$ we deduce that
\[\sum_{i=1}^j \mu_i-\sum_{i=1}^j\lambda_i\le p+q.\]
Combining the above two equations leads to
$\sum_{i=1}^{j} t_i\le p+q$, a contradiction. This completes the proof of the claim. 

Let $\nu=(\mu_1-t_1,\mu_2-t_2,\ldots,\mu_\ell-t_\ell)$. For any $1\le j\le \ell$, if $t_j=0$, then we have $\mu_j-t_j\ge\mu_{j+1}-t_{j+1}$. If $t_j\ge 1$, 
from \eqref{eq-tj} and \eqref{eq-tj-1} it follows that
\[\sum_{i=1}^j(\mu_i-t_i)=\sum_{i=1}^j\lambda_i,\quad \sum_{i=1}^{j-1}(\mu_i-t_i)\le\sum_{i=1}^{j-1}\lambda_i, \quad
\sum_{i=1}^{j+1}(\mu_i-t_i)\le\sum_{i=1}^{j+1}\lambda_i.\]
Hence $\mu_j-t_j\ge\lambda_j\geq \mu_{j+1}-t_{j+1}$. 
This means that $\nu$ is a valid partition. Moreover, according to \eqref{eq-tj-1}, we have 
$\nu \unlhd {\lambda}$. Since $\mathbf{m}\times \mathbf{n}$ is nice, there exists a partition of $\mathbf{m}\times \mathbf{n}$ 
into chains $C_1,C_2\ldots,C_{\ell}$ such that $C_i$ is of size $\nu_i$. Take a partition of $\mathbf{p}\oplus \mathbf{q}$ 
into chains $C'_1,\ldots, C'_{\ell}$ such that $C'_i$ is of size $t_i$. (The existence of such a chain partition is clear since $\mathbf{p}\oplus \mathbf{q}$ is a chain.) Note that for each $1\leq i\leq \ell$ the elements in $C_i \cup C'_i$ form a chain of size $\mu_i$
in $\mathbf{p}\oplus (\mathbf{m}\times \mathbf{n})\oplus \mathbf{q}$. In this way, we obtain a chain partition of $\mathbf{p}\oplus (\mathbf{m}\times \mathbf{n})\oplus \mathbf{q}$ of type $\mu$. This completes the proof.
\end{proof}

The following result can be considered as a generalization 
of Theorem \ref{main-theorem-section4} in some sense. 

\begin{corollary}\label{main-theorem-section4-coro}
For any $p,q\geq 0$, $k\geq 5$ and $n\geq \frac{k+2}{2}$, 
    the distributive lattice $\mathbf{p}\oplus (\mathbf{(n+k)}\times \mathbf{n})\oplus \mathbf{q}$ is not Schur positive.
\end{corollary}

\begin{proof}
Let $\rho$ be the partition $(2n+k-1,2n+k-3,\ldots,k+3,k-3,2,2)$ of $n(n+k)$ as in the proof of Theorem \ref{main-theorem-section4},
and let $\tilde{\rho}$ be the partition of $n(n+k)+p+q$ such that $\tilde{\rho}_1=\rho_1+p+q$ and $\tilde{\rho}_i=\rho_i$ for any $i\geq 2$.
Note that if $[s_{\la}]_{X_{\inc(\PP{(n+k), n})}}\neq 0$, then there exists some $T\in\T_\la$ such that $\SCP_{P,\cont(T)}\neq \emptyset$ by \eqref{formula} and hence $\la_1\leq \cont(T)_1\leq 2n+k-1=\rho_1$ since the longest chain $\PP{(n+k), n}$ is of size $2n+k-1$.
Suppose that 
$$X_{\inc(\PP{(n+k), n})}=\sum_{\substack{\la\vdash n(n+k)\\
\la_1\leq \rho_1}}a_{\la}s_{\la}.$$
It is clear that 
$$X_{\inc(\mathbf{p}\oplus (\mathbf{(n+k)}\times \mathbf{n})\oplus \mathbf{q})}=X_{\inc(\PP{(n+k), n})}\cdot e_{1}^{p+q}.$$
By using Pieri's rule $p+q$ times (see \cite[Theorem 7.15.7]{EC2}), we find that for any $\la\vdash n(n+k)$ with $\la_1\leq \rho_1$ the coefficient $[s_{\tilde{\rho}}]_{(s_{\la} \cdot e_{1}^{p+q})}$ vanishes unless $\la=\rho$. Thus,  
$$[s_{\tilde{\rho}}]_{X_{\inc(\mathbf{p}\oplus (\mathbf{(n+k)}\times \mathbf{n})\oplus \mathbf{q})}}=[s_{\rho}]_{X_{\inc(\PP{(n+k), n})}}<0$$
when $k\geq 5$ and $n\geq \frac{k+2}{2}$.
This completes the proof. 
\end{proof}

\section{Open problems}  \label{section-6}

In this paper, we show that distributive lattices are not nice (thus not Schur positive) in general, giving an answer of Problem \ref{conjecture:Stanley}. However, there are still some open problems about the nice property and the Schur positivity for some special families of distributive lattices, which deserve to be further studied. 

Lonc and Elzobi \cite{LE99} proved that the product of any two chains is nice. It is natural to consider the nice property of the product of more than two chains. Data from computer program show that 
$\mathbf{n_1}\times \mathbf{n_2}\times \cdots\times \mathbf{n_r}$ is nice whenever $n_1n_2\cdots n_r\leq 35$. 
Nevertheless, they are generically not Schur positive by the examples constructed in Section \ref{section-4}. 
Based on computations in SageMath \cite{sage}, we find a bigger number of lattices of the form $\mathbf{m}\times \mathbf{n}$ that are not Schur positive. Together with Theorem \ref{main-theorem-section4}, we have the following conjecture.

\begin{conjecture}\label{conjecture:generalconj}
For any $n\geq 2$, $m\geq 8$ or $n\geq 3$, $m\geq n+5$, 
the distributive lattice ${\inc(\PP{m,n})}$ is not Schur positive. 
\end{conjecture}

To support the above conjecture, let us see some examples. 
When $n=2$, we have verified in SageMath that
\begin{align*}
    [s_{(m+1,m-8,2,2,2,1)}&]_{X_{\inc(\mathbf{m}\times \mathbf{2})}} < 0 \qquad\textnormal{when\ \ \ \ \,} 8\leq m\leq 21,\\
    [s_{(m+1,m-9,2,2,2,2)}&]_{X_{\inc(\mathbf{m}\times \mathbf{2})}} < 0  \qquad\textnormal{when\ \ \ \ \,} 9\leq m\leq 45,\\
    [s_{(m+1,m-10,2,2,2,2,1)}&]_{X_{\inc(\mathbf{m}\times \mathbf{2})}} < 0  \qquad\textnormal{when\ \ \ } 13\leq m\leq 70,\\
    [s_{(m+1,m-11,2,2,2,2,2)}&]_{X_{\inc(\mathbf{m}\times \mathbf{2})}} < 0  \qquad\textnormal{when\ \ \ } 20\leq m\leq 142,\\
    [s_{(m+1,m-12,2,2,2,2,2,1)}&]_{X_{\inc(\mathbf{m}\times \mathbf{2})}} < 0  \qquad\textnormal{when\ \ \ } 25\leq m\leq 227,\\
    [s_{(m+1,m-13,2,2,2,2,2,2)}&]_{X_{\inc(\mathbf{m}\times \mathbf{2})}} < 0  \qquad\textnormal{when\ \ \ } 47\leq m\leq 424,\\
    [s_{(m+1,m-14,2,2,2,2,2,2,1)}&]_{X_{\inc(\mathbf{m}\times \mathbf{2})}} < 0  \qquad\textnormal{when\ \ \ } 64\leq m\leq 667,\\
    [s_{(m+1,m-15,2,2,2,2,2,2,2)}&]_{X_{\inc(\mathbf{m}\times \mathbf{2})}} < 0  \qquad\textnormal{when\ \,} 117\leq m\leq 1199,
\end{align*}
implying that ${\inc(\PP{m,2})}$ is not Schur positive for any $8\leq m\leq 1199$. 
We also have verified that ${\inc(\PP{m,3})}$ is not Schur positive for any $8\leq m\leq 100$ and ${\inc(\PP{m,4})}$ is not Schur positive for any $9\leq m\leq 50$. 

For the product of three chains, we have 
verified that $X_{\inc(\PP{m,2,2})}$ is not Schur positive for $8\leq m\leq 23$.
Our results about the Schur positivity of products of chains tend to imply that, 
for distributive lattices of the form $\mathbf{n_1}\times \mathbf{n_2}\times \cdots\times \mathbf{n_r}$, 
if the chain sizes $n_i$'s have big differences, then the lattice is generically not Schur positive.
Since the boolean algebra $B_r$ is isomorphic to the product of $r$ chains of length $2$, 
Conjecture \ref{conjecture:StanleyBn} suggests that $\mathbf{n_1}\times \mathbf{n_2}\times \cdots\times \mathbf{n_r}$ might be 
Schur positive when $n_i$'s are similar. Motivated by this, we propose the following conjecture. 

\begin{conjecture}\label{conjecture:samechains}
    For any positive integers $n$ and $r$, the product $\mathbf{n}^r=\PP{n,\cdots,n}$ of $r$ chains of length $n$ is Schur positive.
\end{conjecture}

\section*{Acknowledgement}
This work is supported by the Fundamental Research Funds for the Central Universities.
Dun Qiu is supported in part by the National Science Foundation of China (12271023 and 12171034).
Arthur Yang is supported in part by the National Science Foundation of China (11971249 and 12325111).

\end{document}